\documentclass[11pt,twoside]{article}

\usepackage{mathrsfs,amsfonts,amsmath,amssymb}
\usepackage{latexsym}
\newtheorem{satz}{Theorem}

\newtheorem{theorem}[satz]{Theorem}
\newtheorem{lemma}[satz]{Lemma}

\newcommand{\lf}{\left\lfloor}
\newcommand{\rf}{\right\rfloor}

\def\no{\noindent}
\def\sbeq{\subseteq}

\def\Z{\mathbb {Z}}

\def\zn{\Z/N\Z}
\def\e{\varepsilon}

\def\a{\alpha}

\def\C{\mathbb{C}}

\def\h{\widehat}

\def\d{\delta}

\def\({\big (}
\def\){\big )}
\def\g{\gamma}
\def\G{\Gamma}

\def\b{\beta}
\def\le{\leqslant}
\def\ge{\geqslant}
\def\dim{{\rm dim}}
\def\le{\leqslant}
\def\ge{\geqslant}
\def\_phi{\varphi}

\def\m{\times}

\def\D{\Delta}

\def\th{\theta}

\def\th{\theta}

\def\rk{{\rm{ rk}}}
\def\m{\mu}

\def\Span{{\rm Span\,}}

\def\D{\Delta}

\def\La{\Lambda}

\def\hA{\h {1_A}}
\def\hB{\h {1_B}}
\def\hX{\h {1_X}}

\begin{document}
\title{\bf A subexponential upper bound for the van der Waerden numbers $W(3,k)$}

\author{ By\\  \\{\sc Tomasz Schoen}}

\date{}

\maketitle

\begin{abstract} We show an improved upper estimate for the van der Waerden number $W(3,k):$ there is  an absolute  constant $c>0$ such that if  $\{1,\dots,N\}=X\cup Y$ is a partition such that $X$  does not contain any  arithmetic progression of length $3$ and $Y$ does not contain any  arithmetic progression of length $k$ then
$$N\le \exp(O(k^{1-c}))\,.$$ 
\end{abstract}

\section{Introduction}\label{s:intro}

Let $k$ and $l$ be positive integers. The van der Waerden number $W(k,l)$ is the smallest positive integer $N$ such that in any partition
$\{1,\dots, N\}=X\cup Y$ there is an arithmetic progression of length $k$ in $X$ or an arithmetic progression of length $l$ in $Y$. The existence of such numbers was established by van der Waerden \cite{vdW}, however the order of magnitude of $W(k,l)$ is unknown for  $k,l\ge 3.$ Clearly,  $W(k,l)$ is related to Szemer\'edi's theorem on arithmetic progressions \cite{szemeredi} and any effective estimate in this theorem leads to an upper bound on the van der Waerden numbers. Currently best known bounds in the most important 
 diagonal case are
$$(1-o(1))\frac{2^{k-1}}{ek}\le W(k,k)\le {2^{2^{2^{2^{2^{k+9}}}}}}\,.$$  
The upper bound follows from the famous work of Gowers \cite{gowers} and the lower bound was proved by Szab\'o \cite{szabo} using probabilistic argument. 
Furthermore, Berlekamp \cite{berlekamp} showed that if $k-1$ is a prime number then
$$W(k,k)\ge (k-1)2^{k-1}\,.$$ 

 Another very intriguing instance  the problem are numbers $W(3,k)$ as they are related to Roth's theorem \cite{roth} that provides  more efficient estimates for sets avoiding three-terms arithmetic progressions. Let us denote by $r(N)$ the size of the largest progression-free subset of  $\{1,\dots, N\}$. We know that 
\begin{equation}\label{roth}
r(N)\ll \frac{N}{{\log N}^{1-o(1)}}\,,
\end{equation}
see \cite{bloom, bloom-sisask,sanders,schoen}. 
However this bound is not strong enough to imply a subexponential estimate for  $W(3,k)$.

Green \cite{green} proposed a  very clever argument based on  arithmetic properties of sumsets to bound $W(3,k).$ Building on this method and applying results from \cite{crs} it was showed in \cite{cwalina-schoen} that
$$W(3,k)\le \exp(O(k\log k))\,.$$ 
The best known lower bound was obtained by Li and Shu \cite{ls} (see also \cite{blr}), who showed that
$$W(3,k)\gg \Big(\frac{k}{\log k}\Big )^2.$$
The purpose of this paper is to prove a subexponential bound on $W(3,k)$. 

\begin{theorem}\label{main} There are absolute constants $C,c>0$ such that for every $k$ we have
$$W(3,k)\le \exp(Ck^{1-c})\,.$$
\end{theorem}

Our argument is based on the method of \cite{schoen}, which explores in details the structure of a large spectrum.  This method can be partly applied (see Lemma \ref{increment}) in our approach and it deals only with a progression-free  partition class. The second part of the proof exploits the structure of both partition classes  and in this case   the argument of  \cite{schoen} has to be significantly modified.

\section{Notation} \label{s:notation}

The Fourier coefficients of a function $f:\Z/N\Z\to\C$ are defined
by
$$
\h f(r)=\sum_{x=0}^{N-1}f(x)e^{-2\pi ixr/N},
$$
where $r\in\Z/N\Z$. The inversion formula states that
$$f(x)=\frac1N\sum_{x=0}^{N-1}\h f(r)e^{2\pi ixr/N}.$$ We denote by $1_A(x)$  the indicator function of  set $A$.
Thus using the inversion formula and the fact that $\h{(1_A*1_B)}(r)=\hA(r)\hB(r)$ one can express the number of three--term arithmetic progressions (including trivial ones) by
$$\frac1{N}\sum_{r=0}^{N-1}\hA(r)^2\hA(-2r)=|A|\,.$$
 Parseval's identity asserts in particular  that 
$$\sum_{r=0}^{N-1}|\hA(r)|^{2}=|A|N\,.$$
Let  $\th\ge 0$ be a real number,  the $\th-$spectrum of  $A$ is defined by
$$
\D_{\th}(A)=\big \{r\in\Z/N\Z:|\hA(r)|\ge\th|A|\big \}.
$$
If $A$ is specified then we  write $\D_{\th}$ instead of $\D_{\th}(A).$ 

By the span of a finite set $S$ we mean
$$\Span(S)=\Big\{\sum_{s\in S}\e_ss: \e_s\in \{-1,0,1\} \text{ for all } s\in S\Big\}$$
and the dimension of $A$ is defined by
$$\dim(A)=\min\big\{|S|: A\sbeq \Span (S) \big\}\,.$$   
Chang's Spectral Lemma provides an upper bound 
for the dimension of a spectrum.

\begin{lemma} {\rm\cite{chang}}\label{chang}
Let $A\subseteq \Z/N\Z$ be a set of size $|A|=\d N$ and let $\th>0.$ Then
$$\dim(\D_\th(A))\ll \th^{-2}\log(1/\d)\,.$$
\end{lemma}
We are going to use  Bohr sets \cite{bourgain-1/2} to prove the main result.
Let $\G\sbeq\h G$ and   $\g\in(0,\frac{1}{2}]$ then the Bohr set generated by $\G$  with radius $\g$  is 
$$
B(\G,\g)
  =  \big\{ x\in \Z/N\Z: \ \|{tx}/{N}\|\le\g \text { for all } t\in\G\big\}\,,$$
where $\left\Vert x\right\Vert =\min_{y\in\Z}|x-y|$. 
 The rank of $B$ is the size of $\G$ and we denote it by $\rk(B).$
Given $\eta>0$ and a
Bohr set $B=B(\G,\g)$, by $B_{\eta}$ we mean the Bohr set $B(\G,\eta\g).$
We will use two basic properties of Bohr sets concerning its size and regularity, see \cite{tao-vu}.

\begin{lemma} {\rm\cite{bourgain-1/2}}\label{bohr-size}
 For every $\g\in(0,\frac{1}{2}]$ we have
$$
\g^{|\G|}N\le|B(\G,\g)|\le 8^{|\G|+1}|B(\G,\g/2)|\,.
$$
\end{lemma}
We call a Bohr set $B(\G,\g)$ {regular} if for every $\eta$,
where $|\eta|\le 1/(100|\G|)$ we have
$$
(1-100|\G||\eta|)|B|\le|B_{1+\eta}|\le(1+100|\G||\eta|)|B|.
$$
Bourgain \cite{bourgain-1/2} showed that regular Bohr sets are ubiquitous.
\begin{lemma} {\rm\cite{bourgain-1/2}}\label{bohr-regular} For every Bohr set $B(\G,\g),$ there exists
$\g'$ such that $\frac{1}{2}\g\le\g'\le\g$ and $B(\G,\g')$ is
regular.
\end{lemma}

\section{Proof of Theorem \ref{main}}
\label{s:large-coeff}

Our main tool is the next lemma, which can be extracted from \cite{schoen} (see Lemmas 7, 9, 12 and 13). Its proof  makes use of the deep result by Bateman and Katz in \cite{bateman-katz, bateman-katz-nonsmooth} describing  the structure of the large 
spectrum.

\begin{lemma} {\rm\cite{schoen}}\label{increment}
There exists an absolute constant $c>0$ such that   the following holds.
 Let $A\subseteq \Z/N\Z,\, |A|=\d N$ be a set without any non-trivial arithmetic progressions of length three and such that 
\begin{equation}\label{mid}
\sum_{r:\, \d^{1+c} |A|\le |\hA(r)|\le \d^{1/10} |A|}|\hA(r)|^3\ge \frac1{10}\d^{c/5}|A|^3\,.
\end{equation} 
Then  there is a regular Bohr set $B$ with 
$\rk( B)\ll \d^{-1+c}$  and  radius $\Omega( \d^{1-c})$ such that     for some $t$
    $$|(A+t)\cap B|\gg \d^{1-c}|B|.$$
\end{lemma}
		
Furthermore, we  apply Bloom's iterative lemma, that provides a density increment by  a constant factor greater than $1$ for progression-free sets and 	Sanders' lemma on a containment of long arithmetic progressions in dense subsets of regular Bohr sets.
		
\begin{lemma}\label{bloom} {\rm\cite{bloom}} There exists an absolute constant $c_1 > 0$ such that the following
holds. Let $B\sbeq \Z/N\Z$ be a regular Bohr set of rank $d$. Let
$A_1 \sbeq B$ and $A_2 \sbeq B_\e,$ each with relative densities
$\a_i$. Let $\a = \min(c_1, \a_1, \a_2)$ and assume that   $d \le
\exp(c_1(\log^2(1/\a)).$ Suppose that   $B_\e$ is also regular and
$c_1\a/(4d) \le \e \le c_1\a/d.$ Then either
\begin{itemize}
\item[(i)] there is a regular Bohr set $B'$ of rank $\rk(B') \le d +
O(\a^{-1}\log(1/\a))$ and size
$$ |B'| \ge \exp\big(-O(\log^2(1/\a)(d + \a^{-1}\log(1/\a)))\big)|B|$$
such that 
$$|(A_1+t)\cap B'|\gg (1+c_1)\a_1|B'|$$
for some $t\in \zn$;
\item[(ii)] or there are $\Omega(\a_1^2\a_2 |B| |B_\e|)$ three-term arithmetic progressions $x+y=2z$ with $x,y\in A_1, z\in A_2$;
\end{itemize}
\end{lemma}

\begin{lemma} {\rm\cite{sanders}}\label{sanders}
Let $B(\G,\g)\sbeq \zn$  be a regular Bohr set of rank $d$   and let $\e$ be a positive number satisfying  $\e^{-1}\ll \g d^{-1}N^{1/d}$. Suppose
that  $A\sbeq B$ contains at least a proportion $1-\e$ of $B(\G, \g).$ Then 
  $A$ contains an arithmetic
progression of length at least $1/(4\e).$
\end{lemma}

\no {\bf Proof of Theorem \ref{main}.}  Put $M=W(3,k)-1$ and let  $\{1,\dots,M\}=X\cup Y$ be a partition such that $X$ and $Y$ avoid $3$ and $k$-term arithmetic progressions respectively. Clearly, we may assume that $M\ge 100k$ hence
$$|Y|\le M-\lf M/k\rf\le M-M/(2k),$$
as no block of $k$ consecutive numbers is contained in $Y$ and therefore $|X|\ge M/(2k).$
Let $N$ be any prime number satisfying $2M<N\le 4M.$ We embed $\{1,\dots,M\}=X\cup Y$ in $\zn$ in a natural way and observe that $X$ and $Y$ possess the same properties. Put $|X|=\d N$. 
First let us assume that
\begin{equation}\label{mid}
\sum_{r:\, \d^{1+\m} |X|\le |\hX(r)|\le \d^{1/10} |X|}|\hX(r)|^3\ge \d^{\m/5}|X|^3\,,
\end{equation}
Then by Lemma \ref{increment}   there is  $t\in \zn$ and a regular Bohr set $B^0$ with
$\rk( B^0)\ll \d^{-1+c}$  and  radius $\Omega( \d^{1-c})$ such that     
    $$|(X+t)\cap B^0|\gg \d^{1-c}|B|$$
		for some absolute constant $c>0.$
 Writing  $A_0=(X+t)\cap B_0$ we have
$$|X_0\cap B^0|\gg \a|B^0|\,,$$
where 
$$\a\gg \d^{1-c}\,.$$
By Lemma \ref{bohr-size} we have 
$$|B^0|\ge \exp\big (-O(\d^{-1+c}\log(1/\d))\big)N\,.$$
Next, we 
iteratively  apply  Lemma \ref{bloom}. Since after each step the density increases by factor $1+c_1$ it follows that  after  $l\ll \log (1/\a)$ steps  case $(ii)$ of Lemma \ref{bloom} holds.   Let $B^i$ be Bohr sets obtained in the iterative procedure and observe that $\rk(B^i)\ll \a^{-1}\log^2(1/\a)$ for every $i\le l$.  Therefore, there are
$$\Omega( \a^3 |B^k||B^l_\e|)$$
three-term arithmetic progressions in $X,$ where $\e\ge c_1\a/(4\rk(B^l))\gg \a^2\log^2(1/\a).$ 
By Lemma \ref{bloom} and  Lemma \ref{bohr-size}
we have 
$$|B^l|\ge \exp\big(-O(\a^{-1}\log^4(1/\a))\big)N\ge  \exp\big(-O(\d^{-1+c}\log^4(1/\d))\big)N\,,$$
and 
\begin{eqnarray*}
|B^l_\e|&\ge& \exp\big(-O(\a^{-1}\log^3(1/\a))\big)\exp\big(-O(\a^{-1}\log^4(1/\a))\big)N\\
&\ge&  \exp\big(-O(\d^{-1+c}\log^4(1/\d))\big)N\,.
\end{eqnarray*}
Thus, $X$ contains 
 $$\Omega (\d^{3-3c}
\exp\big(-O(\d^{-1+c}\log^4(1/\d))\big)N^2)$$
arithmetic progressions of length three.
Since there are only $|X|$ trivial progressions in $X$ it follows that
$$|X|\gg \d^{3-3c} \exp\big(-O(\d^{-1+c}\log^4(1/\d))\big)N^2\,,$$
so 
$$W(3,k)\ll N\ll \exp\big(O(\d^{-1+c}\log^4(1/\d)\big)\le \exp\big(O(k^{1-c}\log^4k)\big)\,.$$

Next let us assume that \eqref{mid} does not hold. 
By Chang's lemma 
$$d:=\dim(\D_{\d^{1/10}}(X))\ll \d^{-1/5}\log (1/\d)$$
 hence there is a set $\La$ such that $|\La|=d$ and 
$\D_{\d^{1/10}}\sbeq \Span(\La)$. By Lemma \ref{bohr-regular} there is   a regular Bohr set $B=B(\Lambda,\g)$   with radius $\g\gg \d^3$. Let $\b=\frac1{|B|}1_B$ then for every $r\in \D_{\d^{1/10}}$ we have
\begin{equation}\label{beta}
\big|\h \b(r)-1\big|\le \frac1{|B|}\sum_{b\in B} |e^{-2\pi i\lambda b/N}-1|\le \frac{2\pi}{|B|}\sum_{b\in B} \sum_{\lambda\in \La}\|rb/N\|\le 2\pi \d^2\,,
\end{equation}
and similarly $|\h \b(2r)-1| \ll \d^2.$
For $t\in \zn$ put
$$f(t)=\b*1_X(t)\,$$
and note that if for some $t\in [4\g M,(1-4\g)M]$ we have   $f(t)=\frac1{|B|}|X\cap (B+t)|\le \d^{1+c'},$ where  $c'=c/20,$ then since $B+t\sbeq [1,M]$ it follows that
$$|Y\cap (B+t)|\ge 1-\d^{1+c'}\,.$$ 
Therefore, by Lemma \ref{sanders} either $\d^{-1-c'}\gg \g d^{-1}N^{1/d}$ or $Y$ contains an arithmetic progression of length
$\frac14\d^{-1-c'}$. The former inequality implies that
$$k^{1+c'}\gg\g d^{-1}N^{1/d}\gg \d^4N^{O(\d^{1/5}\log^{-1} (1/\d))}\gg k^{-4} N^{O(k^{-1/5}\log^{-1} k)}\,,$$
so 
$$W(3,k)\ll N\ll \exp\big(O(k^{1/5}\log^2 k)\big)\,.$$ If the second alternative holds then 
$$\frac14\d^{-1-c'}<k$$
hence by \eqref{roth}
$$k^{-1/(1+c')}\ll \d\ll (\log N)^{-1+o(1)}$$
so
$$W(3,k)\le N\ll \exp\big(O(k^{\frac1{1+c'}+o(1)})\big)\,.$$
Finally we can assume that for every $t\in [4\g M,(1-4\g)M]$ we have   $f(t)\ge \d^{1+c'}.$ 
Let $T(X)$ denote the number of three-term arithmetic progression in $X$ and let 
$$T(f)=\sum_{x+y=2z}f(x)f(y)f(z)\,.$$ 
Then clearly 
\begin{equation}\label{f}
T(f)\gg \d^{3+3c'}M^2\gg \d^{3+c/6}N^2
\end{equation} 
and we will show that  $T(X)$  does not differ  much from $T(f)$
\begin{eqnarray}\label{t}
\big|T(X)-T(f)\big|&=&\frac1{N}\big|\sum_{r=0}^{N-1} \h {1_X}(r)^2\hX(-2r)-\sum_{r=0}^{N-1} \h f(r)^2\h f(-2r)\big|\nonumber\\
&\le&\frac1{N}\sum_{r=0}^{N-1} |\hX(r)^2\hX(-2r)(1-\h \b(r)^2\h \b (-2r))|\\
&=& S_1+S_2+S_3,\nonumber
\end{eqnarray}
where $S_1,S_2$ and $S_3$ summations of \eqref{t} respectively over $\D_{\d^{1/10}},\D_{\d^{1+c}}\setminus \D_{\d^{1/10}}$ and $\zn\setminus \D_{\d^{1+c}}.$
By \eqref{beta}, \eqref{mid}, Parseval's formula and  H\"older's inequality  we have
$$S_1\ll \d^2 \frac1{N}\sum_{r\in \D_{\d^{1/10}}}|\hX(r)|^3\le \d^3\sum_{r=0}^{N-1} |\h {1_X}(r)|^2= \d^2|X|^2\,,$$
$$S_2\le \d^{1+c/5}|X|^2$$
and 
$$S_3\le \frac2{N} \d^{1+c}|X| \sum_{r=0}^{N-1} |\h {1_X}(r)|^2=\d^{1+c}|X|^2\,.$$
Thus, 
$$\big|T(X)-T(f)\big|\ll \d^{3+c/5}N^2,$$
so by \eqref{f} and the fact that $X$ avoids non-trivial three-term arithmetic progression we have
$$|X|=T(X)\gg \d^{3+c/6}N^2\,,$$
hence 
$$W(3,k)\le N\ll \d^{-2-c/6}\ll k^3$$
which concludes the proof.$\hfill\Box$

{}

\bigskip

\no{Faculty of Mathematics and Computer Science,\\ Adam Mickiewicz
University,\\ Umul\-towska 87, 61-614 Pozna\'n, Poland\\} {\tt
schoen@amu.edu.pl}

\end{document}